\documentclass[11pt]{article}
\usepackage{amsmath,amssymb,amsthm,empheq}
\usepackage{comment}
\usepackage{fancybox}
\usepackage{eurosym}
\usepackage{type1cm}
\usepackage{latexsym}
\usepackage{enumerate}
\usepackage{amsfonts}
\usepackage{graphicx}
\usepackage{tikz}
\usepackage{url}
\usepackage{cite}

\newcommand{\bm}[1]{\mbox{\boldmath $#1$}}
\newcommand{\Prob}{\mbox{\rm Pr}}

\renewcommand{\Vec}[1]{\mbox{\boldmath{$#1$}}}
\newcommand{\Cont}{\mbox{\tt cont}}
\newcommand{\Stop}{\mbox{\tt stop}}

\newtheorem{Theorem}{Theorem}[section]

\title{Dynamic Programming  and Linear Programming \\for Odds Problem
\thanks{
This work was supported by JSPS KAKENHI Grant Numbers 
JP26285045, JP26242027, JP20K04973.  
We thank Katsunori Ano, Naoto Miyoshi, and Akifumi Kira
  for extensive discussions.
}
}

\author{Sachika Kurokawa
\thanks{Graduate School of Engineering, 
   Tokyo Institute of Technology 
}
       \and
        Tomomi Matsui 
\thanks{Graduate School of Engineering, 
   Tokyo Institute of Technology 
}
}

\begin{document}
\maketitle

\begin{abstract}
This paper discusses the odds problem, 
	proposed by 
	Bruss in 2000, 
	and its variants.
A recurrence relation 
	called a dynamic programming (DP) equation 
	is used to find an optimal stopping policy of the odds problem 
	and its variants.
In 2013, Buchbinder, Jain, and Singh proposed 
	a linear programming (LP) formulation 
	for finding an optimal stopping policy 
	of the  classical secretary problem, 
	which is a special case of the odds problem.	
The proposed linear programming problem, 
	which maximizes the probability of a win, 
	differs from the DP equations known 
	for long time periods.	
This paper shows that an ordinary DP equation 
	is a modification of the dual problem 
	of linear programming including the LP formulation 
	proposed by Buchbinder, Jain, and Singh. 
	
 
\end{abstract}

\section{Introduction}
The odds problem is an optimal stopping problem
	proposed by Bruss~\cite{bruss2000sum}, 
	which includes the classical secretary problem as a special case. 
Some variants of the odds problem 
	are discussed in~\cite{		
	ano2010odds,bruss2000selecting,kurushima2016multiple,
	tamaki2010sum,matsui2014note,matsui2016lower,matsui2017compare}.
Although the odds theorem shown in~\cite{bruss2000sum} and~\cite{bruss2003note}
	gives an optimal policy directly, 
	it is well-known that 
	a backward induction method finds an optimal policy 
	for most of the finite optimal stopping problems~\cite{
		gilbert1966recognizing}. 
The backward induction method was used 
	 for the first time by Cayley~\cite{
	cayley1875mathematical,moser1956problem,ferguson1989solved}.
A recurrence relation appearing in the backward induction method
	is called a dynamic programming (DP) equation.
Recently, Buchbinder, Jain, and Singh~\cite{buchbinder2014secretary} 
	proposed a linear programming (LP) formulation
	for the classical secretary problem.
Their LP formulation does not depend on the backward induction method
	and is essentially different from the DP equation.
The purpose of this paper is to clarify 
	the relation between these two types of programming.
	
The classical secretary problem and its variants are discussed 
	as finite-horizon Markov decision processes and 
	solved through backward induction methods 
	(e.g., see Gilbert and Mosteller~\cite{gilbert1966recognizing}, 
	Ross~\cite{ross2014introduction} Section~I.5, 
	Puterman~\cite{puterman2014markov} Section~4.6.4, 
	and Bertsekas~\cite{bertsekas2017dynamic} Section~3.4 of Vol.~I).
Variations of infinite-horizon optimal stopping problems  are 
	discussed as stationary, infinite-horizon Markov decision processes
	(e.g., see~\cite{ross2014introduction} Section~III.2, 
	\cite{puterman2014markov} Section~7.2.8,  and
	\cite{bertsekas2012dynamic} Section~4.4 of Vol.~II).
An LP reformulation is a popular method 
	for stationary, infinite-horizon 
	Markov decision process models~\cite{manne1960linear,kallenberg1994survey}. 
By contrast, an LP formulation is seldom used 
	for  finite-horizon Markov decision process models 
	(and/or transient Markov policies)
	owing to the ease of solving the DP equation.
Theoretically, Derman~\cite{derman1962sequential,derman1965some}
	showed that we can deal with a finite-horizon Markov decision process model
	as a special case of an infinite-horizon Markov decision process models. 
Some recent results on  LP formulations 
	for finite-horizon Markov decision process models
	appear in~\cite{bhattacharya2017linear,mundhenk2000complexity}, 
	for example.
In Section~\ref{DP2LP}, 
	we propose a linear programming problem 
	whose unique optimal solution is the solution of the DP equation.
We give a direct proof of a transformation (from the DP equation 
	to an LP formulation) in a self-contained manner 
	by restricting to the odds problem.

In Section~\ref{FlowF}, 
	we describe the odds problem as a 
	finite-horizon Markov decision process
	and discuss a problem of  finding an optimal policy 
	(stopping rule).
A straightforward formulation, 
	which is independent of the DP equation,
	produces a non-linear programming problem.
We propose a transformation that converts 
	the non-linear programming problem
	into a linear programming problem, called a {\em flow formulation}.
By restricting to the classical secretary problem, 
	our flow formulation gives the LP formulation 
	proposed in~\cite{buchbinder2014secretary}. 
As a main result, we  show that our flow formulation is the dual linear programming 
	of the LP formulation obtained in Section~\ref{DP2LP}.

The remainder of this paper is organized as follows. 
In the next section, 
	we provide detailed descriptions of the odds problem
	and its variants.
Section~\ref{DP2LP} 
	shows that the DP equation gives a unique optimal solution 
	of a certain linear programming problem.
In Section~\ref{FlowF}, 
	we describe a linear programming formulation 
	for determining an optimal policy of  the odds problem.
We also show the duality of 
	a DP equation and our linear programming formulation. 

\section{Odds Problem and its Variations}

Let $X_1, X_2, \ldots, X_n$ denote a sequence 
	of independent Bernoulli random variables. 
If $X_i = 1$, we state that the outcome 
	of random variable $X_i$ is a {\em success}.
Otherwise $(X_i = 0)$, 
	we state that the outcome of $X_i$ is a {\em failure}. 
For each $i \in \{1, 2, \ldots , n\}$, 
	the success probability 
	$p_i=\Prob[X_i=1](0 < p_i <1)$ is given.  
We denote the failure probability 
	by $q_i=\Prob[X_i=0]=1-p_i >0 $ 
	and the {\em odds} by $r_i=p_i/q_i$.
A player observes these random variables sequentially
	one by one and is allowed to select 
	the variable when observing a success.
The odds problem, proposed by Bruss~\cite{bruss2000sum},
	maximizes the probability of
	selecting the last success. 
	
In this paper, we discuss the odds problem 
	and some variations 
	(see~\cite{matsui2017compare} for detailed discussion).  
We assume that a sequence of rewards $(R_1, R_2, \ldots ,R_n)$ 
	is given and a player receives reward $R_i$ 
		when  selecting a success $X_i=1$.
The aim of the player 	is to maximize the expected reward.
The odds problem proposed by Bruss~\cite{bruss2000sum}
	is obtained by setting the reward $R_i$
	to the probability that a success $X_i=1$ is the last success, 
	i.e., $R_i= q_{i+1} \cdot q_{i+2} \cdot \cdots \cdot q_n$
	$(\forall i \in \{1, 2, \ldots ,n\}).$ 
The classical secretary problem is obtained by setting
	$p_i=1/i$ and $R_i=i/n$ 
		$(\forall i \in \{1, 2, \ldots , n\}).$ 
	
Bruss and Paindaveine~\cite{bruss2000selecting} 
	discussed a model in which  a player wants to select
	the final $m$-th success. 
This problem is obtained by setting
\[
	R_i 
	=\left( \prod_{j={i+1}}^n q_j \right) \cdot 
	\left( 
		\sum_{i+1 \leq i_1< \cdots <i_{m-1} \leq n}
		r_{i_1}r_{i_2} \cdots r_{i_{m-1}} 
	\right) \;\;\;
	(\forall i \in \{1, 2, \ldots ,n\}).
\]

Tamaki~\cite{tamaki2010sum} discussed a problem 
	of selecting  any of the last $m$ successes. 
We can express this problem by setting
\[
	R_i 
 	=\left( \prod_{j={i+1}}^n q_j \right) \cdot 
	\left\{ 1+\sum_{h=1}^{m-1}
	\left( 
		\sum_{i+1 \leq i_1< \cdots <i_h \leq n} 
		r_{i_1}r_{i_2} \cdots r_{i_h} 
	\right) 
	\right\} \;\;\;
	(\forall i \in \{1, 2, \ldots ,n\}).
\]

Matsui and Ano~\cite{matsui2017compare} 
	discussed a problem of selecting
	$k$ out of the last $\ell$ successes, 
	where $1 \leq k \leq \ell < N.$
Their model includes the above models as special cases.
The model is obtained by setting 
\[
	R_i 
 	=\left( \prod_{j={i+1}}^n q_j \right) \cdot 
	\left\{ \sum_{h=k-1}^{\ell-1}
			\left(
				\sum_{i+1 \leq i_1< \cdots <i_h \leq n}
				r_{i_1}r_{i_2} \cdots r_{i_h}
			\right) 
	\right\}\;\;\;
	(\forall i \in \{1, 2, \ldots ,n\}).
\]	


\section{Linear Programming Problem for Solving DP Equation}~\label{DP2LP}

This section describes a DP equation 
	for finding an optimal stopping rule of our problem.
We propose a linear programming problem 
	whose unique optimal solution gives the solution 
	to the DP equation.

Let $w_i$ be the expected reward 
	under the condition in which
	a player observes variables $X_1,X_2, \ldots , X_i$; 
	does not select $X_i$; and afterward
	adopts an optimal stopping rule.
Then, $w_0$ denotes the maximum expected reward
	when a player adopts an optimal stopping rule.
The following recurrence relation,
\begin{align}
\begin{aligned} \label{eq:dp}
& w_{i-1} = \max \{q_iw_i + p_i R_i, w_i\} 
	&(i \in \{1,2, \ldots ,n\}), \\
& w_n = 0, \\
\end{aligned}
\end{align}
	is called a DP equation,  
	which calculates an optimal stopping rule 
	through backward induction (see~\cite{gilbert1966recognizing} for example).

Now, we describe a linear programming problem
	that finds a solution to the above DP equation.

\begin{Theorem}
A linear programming problem
\begin{align} \nonumber
\begin{aligned}
\mbox{\rm P: \;\;}
\mbox{\rm min.}~& w_0 \\ 
\mbox{\rm s.t.}~& w_{i-1} \geq q_iw_i+p_i R_i 
											&(\forall i \in \{1,2, \ldots ,n\}), \\
& w_{i-1}\geq w_i &(\forall i\in \{1,2, \ldots ,n\}), \\
& w_n = 0, \\
\end{aligned}
\end{align}
	has a unique optimal solution $\bm{w}^\ast$ 
	that satisfies DP Equation (\ref{eq:dp}).
\end{Theorem}

\noindent
{\bf proof:}
Because problem P has a feasible solution 
	whose objective function is always non-negative, 
	an optimal solution is available.
Let $\bm{w}^\ast$ be an optimal solution of P.
Obviously, $\bm{w}^\ast$ satisfies
\[
w_{i-1}^\ast \geq \max \{q_iw_i^\ast + p_i R_i, w_i^\ast \}
	 \;\;\;(\forall i \in \{1, 2, \ldots ,n\}).
\]
We prove that $\bm{w}^\ast$ satisfies DP Equation (\ref{eq:dp}) 
	through a contradiction.
Assume that there exists an index
	$\hat{i} \in \{1,2, \ldots ,n\}$ that satisfies 
\[
	w_{\hat{i}-1}^\ast >
 \max \{q_{\hat{i}}w_{\hat{i}}^\ast + p_{\hat{i}} R_{\hat{i}}, 
 	w_{\hat{i}}^\ast \}.
\]
We introduce a sufficiently small positive number $\varepsilon>0$ 
	and a solution $(w_0', w_1', \ldots ,w_n')$ defined by 
\[
w_k'=
	\left\{
		\begin{array}{ll}
			w_k^{\ast} & (\forall k \in \{{\hat{i}}, \ldots ,n\} ),\\
			w_k^{\ast}- \varepsilon^{n-k}
										&(\forall k \in \{0,1,\ldots, {\hat{i}}-1\} ).\\
		\end{array}
	\right.
\]
In the following, we show that 
	$(w_0', w_1', \ldots ,w_n')$ is feasible to P.

For each $k \in \{ \hat{i}+1, \hat{i}+2, \ldots,  n\},$
	the definition of ${\bm w'}$ directly implies that 
\begin{align}	\nonumber
\begin{aligned} 
w'_{k-1}=~&w^\ast_{k-1}\geq q_k w^\ast_k +p_k R_k 
					=q_k w'_k +p_k R_k , \;\; \mbox{ and } \\\
w'_{k-1}=~&w^\ast_{k-1} \geq w^\ast_k=w'_k.
\end{aligned}		
\end{align}
When  $k=\hat{i}$, 
	$ w'_{\hat{i}-1}$ satisfies 
\begin{align}	\nonumber
\begin{aligned} 
w'_{\hat{i}-1}=~&w^\ast_{\hat{i}-1}-\varepsilon^{n-\hat{i}+1}
	> q_{\hat{i}} w^\ast_{\hat{i}} +p_{\hat{i}} R_{\hat{i}}
	  =q_{\hat{i}} w'_{\hat{i}} +p_{\hat{i}} R_{\hat{i}}, \;\;  \mbox{ and } \\
w'_{\hat{i}-1}=~&w^\ast_{\hat{i}-1}-\varepsilon^{n-\hat{i}+1} 
  > w^\ast_{\hat{i}}=w'_{\hat{i}},
\end{aligned}		
\end{align}
	as  $\varepsilon$ is a sufficiently small positive number  
	and $n-\hat{i}+1 \geq 1$.
For each \mbox{$k \in \{1, 2,  \ldots ,\hat{i}-1\}$,}  
	the inequality $q_k>0$ implies the following: 
\begin{align}	\nonumber
\begin{aligned} 
w'_{k-1}=~&w^\ast_{k-1}-\varepsilon^{n-(k-1)}
	\geq q_k w^\ast_k +p_k R_k -  \varepsilon^{n-k+1}
  \geq q_k w^\ast_k +p_k R_k -  q_k \varepsilon^{n-k} \\
  =~&q_k w'_k +p_k R_k, \;\;   \mbox{ and } \\
w'_{k-1}=~&w^\ast_{k-1}-\varepsilon^{n-(k-1)}
  \geq w^\ast_k-  \varepsilon^{n-k+1}
  \geq w^\ast_k-  \varepsilon^{n-k}
	=w'_k.
\end{aligned}		
\end{align}
From the above, ${\bm w'}$ is feasible to problem P.
Obviously, we have $w_0'<w_0^\ast$, 
	which contradicts with the optimality of ${\bm w^\ast}$.
	\hfill $\Box$



\section{Flow Formulation and its Duality}\label{FlowF}

\subsection{Flow Formulation}
In this subsection, 
	we describe the odds problem as a finite-horizon 
	Markov decision process 
	and give a straightforward non-linear programming formulation
	of a problem for finding an optimal policy. 
We propose a transformation that
	converts the problem into a linear programming problem.


In this paper, 
	we represent the record of a game 
	using  a sequence of realization values of observed random variables.
For example, 
	the record $(X_1, X_2, X_3, X_4)=(1,0,0,1)$ indicates that 
	the player observed four random variables 
	$X_1, X_2, X_3, X_4$ and 
	selected the success $X_4=1$.
When a player observes all random variables 
	and fails to select one, 
	we define a record of the game though a sequence 
	of realization values  with an additional last component 1, 
	i.e., $(X_1, X_2, \ldots , X_n, 1)$.
Then, the set of all possible records, 
	denoted by~$\Xi$, becomes 
\[
	\Xi =\{(\xi_1, \xi_2 , \ldots , \xi_i) \in \{0, 1\}^i
		\mid i \in \{1,2,\ldots , n+1\}, \xi_i=1\}.
\]

Next, we introduce a finite-horizon Markov decision process
	that formulates our problem.
Let $\varphi$ be a state space defined by 
	$\varphi ={\cal C} \cup {\cal S}$, 
	where ${\cal C}=\{c_0, c_1, \ldots , c_{n+1}\}$ 
	and ${\cal S}=\{s_1, s_2, \ldots , s_n\}$.
The state $c_0$ is the initial state of our process.
The state $c_{n+1}$ is an absorbing state called a {\em terminal state}.
For each state $c_i \; (i\in \{0,1,2, \ldots , n\})$, 
	we define a transition probability
	from $c_i$ to a state $s \in \varphi$ by
\[
	p_{c_i,s}=
	\left\{ \begin{array}{ll}
		p_i & (\mbox{if } s=s_{i+1}),	\\
		q_i & (\mbox{if } s=c_{i+1}), \\
		0   & (\mbox{otherwise}), 
	\end{array} \right.
	(\forall i \in \{0, \ldots n-1\}) \mbox{ and }
	p_{c_n,s}=
	\left\{ \begin{array}{ll}
		1 & (\mbox{if } s=c_{n+1}), \\
		0   & (\mbox{otherwise}).
	\end{array} \right.
\]
For each state in ${\cal S}$, 
	we associate action space ${\cal A}=\{ \Cont, \Stop\}$.
If an action $\Stop$ is selected at state $s_i \in {\cal S}$, 
	the process moves to the terminal state $c_{n+1}$ 
	and generates a reward $R_i$.
Otherwise, an action $\Cont$ is selected and 
the process moves from $s_i$ to $c_{i}$ without a reward.  
A (probabilistic) policy is an $n$-dimensional vector 
	$\pi=(\pi_1, \pi_2, \ldots , \pi_n) \in [0, 1]^n$, 
	where $\pi_i$ denotes a probability that 
	an action $\Cont$ is selected at state $s_i$.
Our goal is to find an optimal policy $\pi^*$ 
	such that the expected 
	reward is maximized.
	
Given a policy $\pi$, 
	the occurrence probability 
	of a record $\xi =(\xi_1, \xi_2, \ldots , \xi_i) \in \Xi$, 
	denoted by $\Prob (\xi \mid \pi)$, satisfies 
\[
	\Prob (\xi \mid \pi)=\alpha_1 \alpha_2 \cdots \alpha_i
\]
	where
\[
	\alpha_j=
	\left\{\begin{array}{ll}
		q_j & (\mbox{if } \xi_j=0), \\
		p_j \pi_j & (\mbox{if } \xi_j=1), 
	\end{array}\right.
	(\forall j \in \{1,2,\ldots ,i-1\})
\]
and
\[
	\alpha_i=
	\left\{\begin{array}{ll}
		p_i (1-\pi_i) & (\mbox{if } i \leq n), \\
		q_n+p_n \pi_n & (\mbox{if } i=n+1).
	\end{array}\right.
\]
For each $i \in \{1, 2,\ldots , n+1\}$, 
	we define 
\[
	\Xi_i=\{ \xi \in \Xi \mid 
		\mbox{the length of $\xi$ is equal to $i$}\}.
\]
Obviously,  the  occurrence probabilities 
$\Prob (\xi \mid \pi)  \; (\xi \in \Xi_i)$ satisfy
\[
	\sum_{\xi \in \Xi_i} \Prob (\xi \mid \pi)
			=\left( \prod_{j=1}^{i-1} (q_j +p_j \pi_j) \right) p_i  (1-\pi_i)
\]
 for each $i \in \{1,2,\ldots , n\}. $ 
The expected reward with respect to $\pi$ is equal to 
\[
	\sum_{i =1}^n \left( 
		R_i \sum_{\xi \in \Xi_i} \Prob (\xi \mid \pi)
		\right).
\]

We can then formulate the problem of 
	maximizing the expected reward as
\begin{align} \nonumber 
\begin{aligned} 
\mbox{Q: } \;\;\mbox{max.}~& \sum_{i=1}^{n} R_i y_i \\
\mbox{s.t.}~& y_i=\left( 
		\prod_{j=1}^{i-1} (q_j +p_j \pi_j) 
	\right) p_i  (1-\pi_i)
			&(\forall i \in \{1,2,\ldots ,n\}), \\
& 0 \leq \pi_i \leq 1
			&(\forall i \in \{1,2, \ldots ,n\} ). \\
\end{aligned}
\end{align}

\noindent
In the following, we transform the above problem 
	into a linear programming problem.
The following theorem provides the key idea of our transformation.

\begin{Theorem} \label{theo:NLP2LP}
Let $\Vec{y}=(y_1, y_2, \ldots , y_n)$ be an $n$-dimensional 
	real vector.
The following two conditions are equivalent.

\begin{description}
\item[(c1)] There exists $\pi=(\pi_1, \pi_2, \ldots , \pi_n)$
	such that $(\Vec{y}, \pi)$ is feasible for problem Q. 
\item[(c2)] There exists $\Vec{z}=(z_0, z_1, \ldots  , z_n)$ 
	such that $(\Vec{y}, \Vec{z})$ satisfies 
	the following linear inequality system:
\begin{align} \label{ineq-Th}
\begin{aligned} 
& y_i \leq p_i z_{i-1} &(\forall i \in \{1,2, \ldots ,n\}), \\
& y_i + z_i = z_{i-1} &(\forall i \in \{1,2, \ldots ,n\}), \\
& z_0 = 1, \\
& y_i \geq 0 &(\forall i \in \{1,2, \ldots ,n\}). \\
\end{aligned}
\end{align}
\end{description}
\end{Theorem}

\noindent
\underline{
Proof of 
(c1) $\rightarrow$ (c2).} 
Let  $(\Vec{y}, \pi)$ be a feasible solution to problem Q. 
Obviously, we have $\Vec{y}\geq \Vec{0}$.
Apply $z_0=1$  and 
	$z_i=\prod_{j=1}^i (q_j+p_j \pi_j)$
	 $(\forall i\in \{1, 2, \ldots , n\})$.
We then obtain for each $i \in \{1, 2, \ldots , n\}, $ 
\begin{eqnarray*}
\lefteqn{
	z_{i-1}-z_i = 
	\prod_{j=1}^{i-1} (q_j+p_j \pi_j)
	- \prod_{j=1}^{i} (q_j+p_j \pi_j)} \\
	&=& \left( \prod_{j=1}^{i-1} (q_j+p_j \pi_j) \right)
		\left( 1- (q_i + p_i \pi_i) \right) 
	= \left( \prod_{j=1}^{i-1} (q_j+p_j \pi_j) \right)
		 p_i(1- \pi_i) 
	=y_i
\end{eqnarray*}
and 
\begin{eqnarray*}
	p_i z_{i-1}-y_i &=&
	p_i \prod_{j=1}^{i-1} (q_j+p_j \pi_j)
	- \left( \prod_{j=1}^{i-1} (q_j+p_j \pi_j) \right)
		 p_i(1- \pi_i)  \\
	&=& \left( \prod_{j=1}^{i-1} (q_j+p_j \pi_j) \right)
		( p_i - p_i (1-\pi_i))
	= \left( \prod_{j=1}^{i-1} (q_j+p_j \pi_j) \right) p_i \pi_i \geq 0.
\end{eqnarray*}
Thus, $(\Vec{y}, \Vec{z})$ satisfies 
	the inequality system~(\ref{ineq-Th}).  

\smallskip 
\noindent
\underline{
Proof of 
(c2) $\rightarrow$ (c1).} 
Assume that  $(\Vec{y}, \Vec{z})$ satisfies
	the inequality system~(\ref{ineq-Th}).
First, we show $z_i>0$ through the induction on $i$. 
Clearly, $z_0=1 >0$, and  $z_{i-1}>0$ implies that
	$z_i=z_{i-1}-y_i \geq z_{i-1} -p_i z_{i-1}=z_{i-1}(1-p_i) >0. $  
With  $p_i >0$, we define
\begin{equation} \label{LP2policy}
	\pi_i= 1- y_i/(p_i z_{i-1}) 	\;\;\; (\forall i \in \{1, 2, \ldots , n\}).
\end{equation}
We then have the following inequalities:
\[
	1 \geq \pi_i=1- \frac{y_i}{p_i z_{i-1}}
	= \frac{p_i z_{i-1}-y_i}{p_i z_{i-1} } \geq 0 	
	\;\;\; (\forall i \in \{1, 2, \ldots , n\}).
\]
It is easy to see that  
	for each $i \in \{1, 2, \ldots ,n\}$,  
\begin{eqnarray*}
	z_i&=&z_{i-1} -y_i = (q_{i} +p_{i})z_{i-1}-y_i 
				=q_{i}z_{i-1} +(p_{i}z_{i-1}-y_i)  \\
		&=& q_{i}z_{i-1} +p_iz_{i-1}\pi_i 
			=(q_i+p_i \pi_i)z_{i-1}  \\
		&=&(q_i+p_i \pi_i)(q_{i-1}+p_{i-1}\pi_{i-1})\cdots (q_1+p_1 \pi_1)z_0 
		=\prod_{j=1}^{i} (q_j + p_j \pi_j), 
\end{eqnarray*}
	and thus, 
\begin{eqnarray*}
y_i &=& p_i  z_{i-1} - p_i \pi_i z_{i-1} = z_{i-1} p_i   (1-\pi_i)
	 = \left( \prod_{j=1}^{i-1} (q_j + p_j \pi_j) \right) p_i(1-\pi_i). 
\end{eqnarray*}
 From the above, $(\Vec{y}, \pi)$ is feasible to problem Q. 
 \hfill $\Box$
 
\medskip
 
By employing the above theorem, 
	we  can transform problem Q 
	into the following linear programming problem:
\begin{align}
\mbox{FF: } \;\;\mbox{max.}~& \sum_{i=1}^{n} R_i y_i \nonumber \\
\mbox{s.t.}~& \frac{y_i}{p_i} -  z_{i-1} \leq 0 
			&(\forall i \in \{1,2,\ldots ,n\}), \nonumber \\
& y_i + z_i - z_{i-1}=0 
			&(\forall i \in \{1,2, \ldots ,n\} ),   \label{flow-con-law} \\ 
& z_0 = 1,  \label{source-flow} \\
& y_i \geq 0 
		&(\forall i \in \{1,2, \ldots ,n\}).  \nonumber
\end{align}
We call the above a {\em flow formulation}.
When we have an optimal solution to FF, 
	Equation~(\ref{LP2policy}) provides an optimal policy $\pi$, 
	which is optimal to problem Q.

In the following, we interpret the FF problem
	as a flow problem on a digraph.
Let $G=(V, E)$ be a digraph
	with a vertex set 
	$V=\{s_0, s_1, \ldots ,s_n \} \cup \{t_1, t_2, \ldots , t_n\}$ 
	and a directed edge set
	$E=\{(s_{i-1}, s_i) \mid i \in \{1, 2, \ldots , n\} \}
		\cup \{(s_{i-1}, t_i) \mid    i \in \{1, 2, \ldots , n\} \}$. 
For each $i \in \{1, 2, \ldots , n\}$,  
	we associate variables  $z_i$ and $y_i$
	with edges $(s_{i-1}, s_i)$  and $(s_{i-1}, t_i)$, respectively.
If we regard these variables 
	as a flow on a corresponding directed edge, 
	constraint~(\ref{flow-con-law})  represents 
	a flow conservation law at vertex $s_{i-1}$,
	 and  
	constraint~(\ref{source-flow}) states that
	a flow of volume 1 emanates from vertex $s_0$.
Because $y_i=\sum_{\xi \in \Xi_i} \Prob (\xi \mid \pi)$, 
	vertex $t_i$ corresponds to an event in which
	a player selects a success $X_i=1$.

\usetikzlibrary{shapes,arrows,positioning,automata}
\begin{figure}[htb]
\begin{center}
\begin{tikzpicture}[->,>=stealth',shorten >=2pt, line width=1pt, 
                 node distance=1.0cm, style ={minimum size=6mm},	
										initial text={$(z_0=1)$}]
\tikzstyle{every node}=[font=\large]
\node [circle, draw, initial] (s0) {$s_0$};
\node [circle, draw] (s1) [right=of s0] {$s_1$};
\node [circle, draw] (s2) [right=of s1] {$s_2$};
\node [circle, draw] (s3) [right=of s2] {$s_3$};
\node [circle, draw] (s4) [right=of s3] {$s_4$};
\node [circle, draw] (s5) [right=of s4] {$s_5$};
\node [circle, draw] (t1) [below=of s1] {$t_1$};
\node [circle, draw] (t2) [below=of s2] {$t_2$};
\node [circle, draw] (t3) [below=of s3] {$t_3$};
\node [circle, draw] (t4) [below=of s4] {$t_4$};
\node [circle, draw] (t5) [below=of s5] {$t_5$};

\node [draw,color=white,text=black, font=\small](stop1)[below= 0cm of t1] {select $X_1$};
\node [draw,color=white,text=black, font=\small](stop2)[below= 0cm of t2] {select $X_2$};
\node [draw,color=white,text=black, font=\small](stop3)[below= 0cm of t3] {select $X_3$};
\node [draw,color=white,text=black, font=\small](stop4)[below= 0cm of t4] {select $X_4$};
\node [draw,color=white,text=black, font=\small](stop5)[below= 0cm of t5] {select $X_5$};

\path[->] (s0) edge node [above] {$z_1$} (s1)
						(s0) edge node [above] {$y_1$} (t1)
						(s1) edge node [above] {$z_2$} (s2)
						(s1) edge node [above] {$y_2$} (t2)
						(s2) edge node [above] {$z_3$} (s3)
						(s2) edge node [above] {$y_3$} (t3)
						(s3) edge node [above] {$z_4$} (s4)
						(s3) edge node [above] {$y_4$} (t4)
						(s4) edge node [above] {$z_5$} (s5)
						(s4) edge node [above] {$y_5$} (t5);
\end{tikzpicture}
\caption{Digraph $G$, where  $n=5.$}
\end{center}
\end{figure}
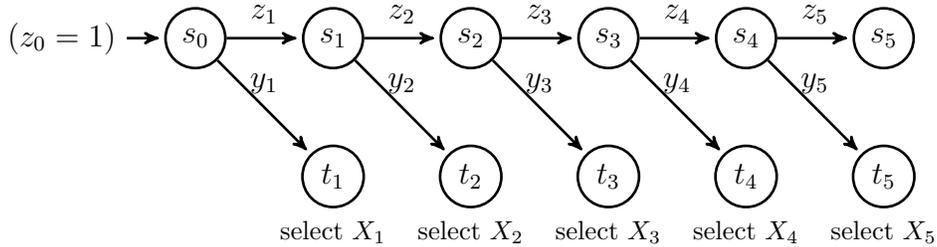

\subsection{Flow Formulation for Classical Secretary Problem}

Buchbinder, Jain, and Singh~\cite{buchbinder2014secretary}
	proposed the following LP formulation
	for the classical secretary problem:  
\begin{align}
\begin{aligned} 
\mbox{max.}~&  \sum_{i=1}^{n} (i/n)  y_i \\
\mbox{s.t.}~& i \cdot y_i \leq  1-\sum_{k=1}^{i-1} y_k  
        &(\forall i\in \{1,2, \ldots ,n\}), \\
& y_i \geq 0 
					&(\forall i \in \{1,2, \ldots ,n\}). \\
\end{aligned}
\end{align}
Their formulation is obtained from the FF
	by eliminating the variables $\{z_1, z_2, \ldots z_n\}$
	through substitutions  $z_i=1-(y_1+y_2+\cdots +y_i)$
	and setting 
	$p_i=1/i$ and $R_i=i/n$ 
		$(\forall i \in \{1, 2, \ldots , n\}).$

\subsection{Dual of Flow Formulation}

The linear programming duality of FF is 
\begin{align}  
\mbox{P1: } \;\;\mbox{min.}~& w_0 \nonumber \\
\mbox{s.t.}~& w_i+\frac{\alpha_i}{p_i} \geq R_i 
			&(\forall i \in \{1,2,\ldots ,n\}), \label{P1-1} \\
& w_{i-1}-w_i-\alpha _i = 0 
			&(\forall i \in \{1,2,\ldots ,n\}),  \nonumber \\
& w_n = 0, \nonumber \\
& \alpha _i \geq 0
			&(\forall i \in \{1,2,\ldots ,n\}).  \label{P1-2} 
\end{align}
We eliminate variables $\{\alpha_1, \alpha_2, \ldots , \alpha_n\}$ 
	by $\alpha_i=w_{i-1}-w_i$.
Constraint~(\ref{P1-1}) is transformed as follows:  
\begin{align}  \nonumber
\begin{aligned}
 w_i+\frac{\alpha_i}{p_i} &\geq &&R_i, \\ 
 p_i w_i+(w_{i-1}-w_i) &\geq &&p_i R_i, \\
 w_{i-1} & \geq && q_i w_i + p_i R_i. \\
\end{aligned}
\end{align}
Non-negativity constraints $\alpha_i \geq 0$,
	as shown in~(\ref{P1-2}), 
	implies the following: 
\[
    w_{i-1}\geq w_i \;\;\; (\forall i \in \{1,2, \ldots ,n\}).
\]
Thus, the above procedure coverts problem P1 into P.

\section{Conclusion}
This paper described a linear programming problem 
	whose unique optimal solution is the solution to the DP equation
	for the odds problem. 
We proposed a flow formulation that maximizes the expected reward
	of the odds problem based on a finite-horizon Markov decision process.
Furthermore, we showed the duality of these two formulations. 




\bibliographystyle{elsarticle-num} 

\bibliography{Odds3}

\begin{thebibliography}{10}
\expandafter\ifx\csname url\endcsname\relax
  \def\url#1{\texttt{#1}}\fi
\expandafter\ifx\csname urlprefix\endcsname\relax\def\urlprefix{URL }\fi
\expandafter\ifx\csname href\endcsname\relax
  \def\href#1#2{#2} \def\path#1{#1}\fi

\bibitem{bruss2000sum}
F.~T. Bruss, Sum the odds to one and stop, Annals of Probability 28~(3) (2000)
  1384--1391.

\bibitem{ano2010odds}
K.~Ano, H.~Kakinuma, N.~Miyoshi, Odds theorem with multiple selection chances,
  Journal of Applied Probability 47~(4) (2010) 1093--1104.

\bibitem{bruss2000selecting}
F.~T. Bruss, D.~Paindaveine, Selecting a sequence of last successes in
  independent trials, Journal of Applied Probability 37~(2) (2000) 389--399.

\bibitem{kurushima2016multiple}
A.~Kurushima, K.~Ano, Multiple stopping odds problem in \mbox{Bernoulli} trials
  with random number of observations, Mathematica Applicanda 44~(1) (2016)
  209--220.

\bibitem{tamaki2010sum}
M.~Tamaki, Sum the multiplicative odds to one and stop, Journal of Applied
  Probability 47~(3) (2010) 761--777.

\bibitem{matsui2014note}
T.~Matsui, K.~Ano, A note on a lower bound for the multiplicative odds theorem
  of optimal stopping, Journal of Applied Probability 51~(3) (2014) 885--889.

\bibitem{matsui2016lower}
T.~Matsui, K.~Ano, Lower bounds for \mbox{Bruss'} odds problem with multiple
  stoppings, Mathematics of Operations Research 41~(2) (2016) 700--714.

\bibitem{matsui2017compare}
T.~Matsui, K.~Ano, Compare the ratio of symmetric polynomials of odds to one
  and stop, Journal of Applied Probability 54~(1) (2017) 12.

\bibitem{bruss2003note}
F.~T. Bruss, A note on bounds for the odds theorem of optimal stopping, Annals
  of Probability 31~(4) (2003) 1859--1961.

\bibitem{gilbert1966recognizing}
J.~P. Gilbert, F.~Mosteller, Recognizing the maximum of a sequence, Journal of
  the American Statistical Association 61~(313) (1966) 35--73.

\bibitem{cayley1875mathematical}
A.~Cayley, Mathematical questions with their solutions, The Educational Times
  23 (1875) 18--19.

\bibitem{moser1956problem}
L.~Moser, On a problem of \mbox{Cayley}, Scripta Math 22 (1956) 289--292.

\bibitem{ferguson1989solved}
T.~S. Ferguson, Who solved the secretary problem?, Statistical Science 4~(3)
  (1989) 282--296.

\bibitem{buchbinder2014secretary}
N.~Buchbinder, K.~Jain, M.~Singh, Secretary problems via linear programming,
  Mathematics of Operations Research 39~(1) (2014) 190--206.

\bibitem{ross2014introduction}
S.~M. Ross, Introduction to stochastic dynamic programming, Academic Press,
  2014.

\bibitem{puterman2014markov}
M.~L. Puterman, Markov decision processes: discrete stochastic dynamic
  programming, John Wiley \& Sons, 2014.

\bibitem{bertsekas2017dynamic}
D.~P. Bertsekas, Dynamic programming and optimal control 4th edition, volume I,
  Athena Scientific, 2017.

\bibitem{bertsekas2012dynamic}
D.~P. Bertsekas, Dynamic programming and optimal control 4th edition, volume
  II, Athena Scientific, 2012.

\bibitem{manne1960linear}
A.~S. Manne, Linear programming and sequential decisions, Management Science
  6~(3) (1960) 259--267.

\bibitem{kallenberg1994survey}
L.~C.~M. Kallenberg, Survey of linear programming for standard and nonstandard
  \mbox{Markovian} control problems. \mbox{Part I: Theory}, Zeitschrift f{\"u}r
  Operations Research 40~(1) (1994) 1--42.

\bibitem{derman1962sequential}
C.~Derman, On sequential decisions and \mbox{Markov} chains, Management Science
  9~(1) (1962) 16--24.

\bibitem{derman1965some}
C.~Derman, M.~Klein, Some remarks on finite horizon \mbox{Markovian} decision
  models, Operations Research 13~(2) (1965) 272--278.

\bibitem{bhattacharya2017linear}
A.~Bhattacharya, J.~P. Kharoufeh, Linear programming formulation for
  non-stationary, finite-horizon \mbox{Markov} decision process models,
  Operations Research Letters 45~(6) (2017) 570--574.

\bibitem{mundhenk2000complexity}
M.~Mundhenk, J.~Goldsmith, C.~Lusena, E.~Allender, Complexity of finite-horizon
  markov decision process problems, Journal of the ACM (JACM) 47~(4) (2000)
  681--720.

\end{thebibliography}





\end{document}